\documentclass[12pt, reqno]{amsart}
\usepackage{amsmath, amstext, amsbsy, amssymb, amscd}

\newtheorem{theorem}{Theorem}[section]
\newtheorem{lemma}[theorem]{Lemma}

\theoremstyle{definition}

\theoremstyle{remark}
\newtheorem{remark}[theorem]{Remark}

\numberwithin{equation}{section}

\newcommand{\bo}{\mathcal B_0}
\newcommand{\binf}{b_\infty}

\newcommand{\dinf}{d_\infty}

\newcommand{\fock}{\mathcal F}
\newcommand{\fockh}{\mathbb F}

\newcommand{\hf}{\frac12}
\newcommand{\hZ}{ \hf +\mathbb Z }
\newcommand{\la}{\lambda}
\newcommand{\OSP}{{\mathcal OSP}}

\newcommand{\SP}{\mathcal SP}
\newcommand{\Tr}{ {\rm Tr} }
\newcommand{\vac}{|0\rangle}
\newcommand{\winfty}{\mathcal W_{1+\infty}}
\newcommand{\Z}{ \mathbb Z }

{\vskip-\lastskip\medskip
  \noindent
  {\em #1.}\enspace
  }%
{\qed\par\medskip
  }

\begin{document}
\title[Correlation functions and twisted Fock spaces]
{Correlation functions of strict partitions and twisted Fock spaces}

\author[Weiqiang Wang]{Weiqiang Wang}
\address{Department of Mathematics, University of Virginia,
Charlottesville, VA 22904} \email{ww9c@virginia.edu}
\thanks{Partially supported by an NSF grant}

\begin{abstract}
Using  twisted Fock spaces, we formulate and study two twisted
versions of the $n$-point correlation functions of Bloch-Okounkov
\cite{BO}, and then identify them with $q$-expectation values of
certain functions on the set of (odd) strict partitions. We find
closed formulas for the $1$-point functions in both cases in terms
of Jacobi $\theta$-functions. These correlation functions afford
several distinct interpretations.
\end{abstract}

\maketitle
\date{}

\section*{Introduction}

Bloch-Okounkov \cite{BO} studied certain correlation functions on
the infinite wedge representation (i.e. the Fock space of a pair
of free fermions), which is closely related to the correlation
functions of vertex operators \cite{Zhu}. Recently it has
gradually become clear that this correlation function (and its
variant) affords several distinct interpretations:

\begin{enumerate}
\item  It can be regarded as a certain character of
representations of the $\winfty$ algebra of level one.

\item It can be regarded as the $q$-expectation value of a certain function on
the set of partitions.

\item It can be interpreted as a generating function of the
stationary Gromov-Witten invariants of an elliptic curve.

\item A variant of this correlation function can be interpreted as
a generating function of equivariant intersection numbers on
Hilbert schemes of points on the affine plane.

\item It can be interpreted as a generating function of certain
structure constants of the class algebras of the symmetric groups.
\end{enumerate}

The items (1) and (2) were points of view taken in \cite{BO} (also
see \cite{Ok}). More on representation theory of $\winfty$ (which
is a central extension of the Lie algebra of differential
operators on the circle) of positive integral levels using Fock
spaces can be found in \cite{FKRW, AFMO}. Item (3) is due to the
recent remarkable work of Okounkov-Pandharipande \cite{OP}. Item
(4) is due to Li, Qin, and the author \cite{LQW}, and (5) is
essentially equivalent to (4) based on the results of Vasserot,
Lascoux-Thibon and ours (cf. \cite{LQW}; also see \cite{LT, Wa}).
We refer to \cite{LT, Wa} for a closely related study of the class
algebras of the symmetric groups and more generally of wreath
products and their connections to the $\winfty$ algebra.

The purpose of this paper is to present (two) twisted versions of
the above correlation functions. This is also partly motivated by
the representation theory of distinguished Lie subalgebras
$\widehat{\mathcal D}^\pm$ of $\winfty$ \cite{KWY} via the twisted
fermionic Fock spaces of Neveu-Schwarz type and of Ramond type,
which are intimately related to Lie subalgebras of
$\widehat{gl}_\infty$ \cite{DJKM}. One can interpret our
correlation functions in a way analogous to (1), (2), and (5)
above, although the connections with the counterpart of (5) will
not be presented in this paper. We believe that there will be a
geometric interpretation of our correlation functions.

The first construction uses the twisted Fock space $\fock$ of a
fermion with integer indices (i.e. the Ramond sector). The
$n$-point correlation functions of certain distinguished operators
(quadratic in fermion generators) in $\fock$ are further
identified as the $q$-expectation value of certain functions on
the set of {\em strict} partitions. The second construction uses
the twisted Fock space $\fockh$ of a fermion with half-odd-integer
indices (i.e. the Neveu-Schwarz sector). The $n$-point correlation
functions of certain operators in $\fockh$ are then interpreted as
the $q$-expectation value of functions on the set of {\em odd
strict} partitions. Using combinatorial methods including an
identity of Euler on partitions, we derive clean closed formulas
of the the $1$-point functions in terms of Jacobi
$\theta$-functions in both cases. From the combinatorial
viewpoint, the representation theory serves as a motivation to
suggest ``right" combinatorial questions (which are justified by
nice answers). Our results should have implications on random
partition theory (compare \cite{Ok}).

In both cases, we obtain the $q$-difference equations for the
$n$-point correlation functions and our formulation also indicates
clearly the nature of poles of the $n$-point correlation
functions. These results have their counterparts in \cite{BO} and
our method is analogous to \cite{Ok}. Indeed, such properties play
a key role in \cite{BO} to establish closed formulas for their
$n$-point correlation functions. We expect these properties will
also play an important role in eventually finding closed formulas
for the $n$-point functions studied in the paper. Our correlation
function also admits a natural {\em super} counterpart (cf.
Remark~\ref{rem:super}).

The paper is organized as follows. In Sect.~\ref{sec:fock1}, we
study a distinguished operator on $\fock$ and formulate the
correlation functions. In Sect.~\ref{sec:strict}, we use a
combinatorial method to derive the $1$-point function. In
Sect.~\ref{sec:property} we establish a $q$-difference equation
for the $n$-point functions. Sect.~\ref{sec:fock2} is a
counterpart of Sect.~\ref{sec:fock1}, where the setup is switched
to the Fock space $\fockh$. In Sect.~\ref{sec:oddstr}, we
formulate and establish the counterparts in $\fockh$ of
Sections~\ref{sec:strict} and \ref{sec:property}. We end in
Sect.~\ref{sec:discuss} with discussions.

\vspace{.1in}

\noindent {\bf Notations.} We denote by $\Z$  the set of integers,
$\Z_+$ the set of non-negative integers, and $\hZ$ the set of
half-odd-integers.

\vspace{.1in}

\noindent {\bf Acknowledgment.} After we finished this work, a
very interesting paper of Milas \cite{Mi} appeared in
arXiv:math.QA/0303154, whose results and approach are independent
of and complementary to ours. He systematically used vertex
operator techniques to study closely related correlation functions
and established precise relations between \cite{BO} and
\cite{Zhu}, while the methods in this paper are mainly
combinatorial and Lie theoretic.

\section{Correlation functions on the Fock space of Ramond type}
\label{sec:fock1}

Introduce a fermionic field of Ramond type
$$\phi (z) =\sum_{n \in \Z} \phi_n z^{-n}$$
with the following commutation relations:
\begin{eqnarray*}
 [ \phi_m, \phi_n ]_+ = \delta_{m,-n}, \quad m, n \in \Z.
\end{eqnarray*}
Note that $\phi_0^2 = 1/2$, and $\phi_0$ anticommutes with all
$\phi_n, n\neq 0$.

Denote by $\fock$ the Fock space of $\phi (z)$ with the highest
weight vector $\vac$ annihilated by $\phi_n$, $n >0$. Recall a
partition $\la =(\la_1, \la_2, \ldots \la_\ell)$ is {\em strict}
if $\la_1 > \la_2 > \ldots > \la_\ell>0$. The number $l$ is often
referred to as the length of $\la$ and denoted by $\ell (\la)$.
Introduce the notations
\begin{eqnarray*}
\phi_\la &=& \phi_{-\la_1} \phi_{-\la_2} \ldots
\phi_{-\la_\ell}  \vac\\
\phi_\la \phi_0 &=& \phi_{-\la_1} \phi_{-\la_2} \ldots
\phi_{-\la_\ell} \phi_0 \vac.
\end{eqnarray*}
Then $\fock$ has a linear basis given by $\phi_\la $ and $\phi_\la
\phi_0$, where $\la$ runs over the set $\SP$ of all strict
partitions (including the empty partition). The following graded
operator $L_0$ (which can be identified with the zero-mode of a
Virasoro algebra) defines a natural $\Z_+$-grading on $\fock$:
$$L_0 \phi_\la
= |\la|  \cdot \phi_\la, \quad
L_0 \phi_\la \phi_0= |\la|  \cdot \phi_\la\phi_0$$
for all $\la \in \SP$. A subalgebra $\binf$ of
$\widehat{gl}_\infty$ acts on $\fock$ as follows \cite{DJKM} (also
cf. \cite{KWY}). The algebra $\binf$ is spanned by $E_{i,j}
-E_{-j,-i}$, where $i,j\in\Z$, and a central element $C$, and the
action is given by
\begin{eqnarray} \label{eq:action}
E_{i,j} -E_{-j,-i} =\; :\phi_{-i}\phi_j:
\end{eqnarray}
and $C$ acts as the identity operator $I$. Recall that the normal
ordered product $ :\phi_{-i}\phi_j:$ equals  $\phi_{-i}\phi_j$
unless $i=j<0$; if $i=j<0$, it equals $\phi_{-i}\phi_j -1.$

It is convenient to introduce the generating field
$$ B(z,w) =\sum_{i,j \in \Z} (E_{i,j} -E_{-j,-i}) z^i w^{-j}.$$
Then the representation of $\binf$ on $\fock$ can be concisely
formulated by
$$B(z,w) = :\phi(z) \phi(w):. $$

The Fock space $\fock$ decomposes into a direct sum of two
irreducible representations of $\binf$:
$$\fock =\fock_0 \bigoplus \fock_1.$$
Here $\fock_0$ (resp. $\fock_1$) are spanned by $\phi_{\la}$ with
$\ell (\la)$ even (resp. odd) and $\phi_{\la} \phi_0$ with $\ell
(\la)$ odd (resp. even), where $\la \in \SP$.

Introduce the following operator in $\binf$ (where $t$ is some
variable):
\begin{eqnarray}  \label{eq:keyop}
  \mathcal B_0 (t) = \sum_{k >0} \left( t^k -t^{-k}
   \right)
  (E_{k,k} -E_{-k,-k}) + \frac{t+1}{2(t-1)} \text{I}
\end{eqnarray}
Note that $\frac{t+1}{2(t-1)}$ can be rewritten in a more
symmetric form as $\frac{t^{1/2} +t^{-1/2}}{2(t^{1/2} -t^{-1/2})
}.$

The operator $\mathcal B_0(t)$ can be formally understood as
$$\mathcal B_0(t)= \sum_{k \in \Z} t^k \phi_{-k}\phi_k$$
which is the zero-mode of $\phi(tz) \phi(z)$ without normal
ordering. Note that  for $|t|>1$, we have
\begin{eqnarray*}
\mathcal B_0 (t) &= & :\mathcal B_0 (t): +
 \left ( 1/2 +\sum_{k<0} t^k \right ) \text{I} \\
  &= & :\mathcal B_0 (t): + \frac{t+1}{2(t-1)} \;  \text{I}
\end{eqnarray*}
where
$$:\mathcal B_0 (t):
 =\sum_{k >0} \left( t^k -t^{-k}
   \right):\phi_{-k}\phi_k:.$$

We introduce the {\em $n$-point correlation functions}
\begin{eqnarray*}
 R(t_1, \ldots, t_n) &=&\Tr \mid_{\fock_0} q^{L_0} \mathcal B_0
(t_1) \cdots \mathcal B_0 (t_n).
\end{eqnarray*}

Similarly we define
\begin{eqnarray*}
 :R:(t_1, \ldots, t_n) &=& \Tr \mid_{\fock_0} q^{L_0} :\mathcal B_0
(t_1): \cdots :\mathcal B_0 (t_n):.
\end{eqnarray*}
It is clear that a multiple of $2$ shows up if we have used $\Tr
\mid_{\fock}$ instead of $\Tr \mid_{\fock_0}$ above.

Apparently, there is a simple identity to express $R(t_1, \ldots,
t_n)$ in terms of $:R:(t_1, \ldots, t_n)$, and vice versa.
Further, for $\la\in \SP$,
\begin{eqnarray}  \label{eq:bridge}
:\mathcal B_0 (t): (\phi_\la ) &=& \sum_{k \ge 1} (t^{\la_k}
-t^{-\la_k}) \cdot \phi_\la   \nonumber \\
:\mathcal B_0 (t): (\phi_\la \phi_0) &=& \sum_{k \ge 1} (t^{\la_k}
-t^{-\la_k}) \cdot \phi_\la \phi_0.
\end{eqnarray}

\section{The correlation function via strict partitions}
\label{sec:strict}

Denote by $\SP_n$ the set of strict partitions of $n$, and thus
$\SP = \cup_{n=0}^\infty \SP_n.$ Denote by
$$(a;q)_r =(1-a)(1-aq)\cdots (1-aq^{r-1}), \quad
(a;q)_\infty =\prod_{i=0}^\infty (1-a q^i). $$ It is well known
that
$$\sum_{\la \in \SP} q^{|\la|} =
(-q;q)_\infty = (q;q^2)_\infty^{-1}. $$
Given a function $f$ on the set ${\mathcal SP}$, we denote by
$$\langle f \rangle^{\SP}_q =(-q;q)_\infty^{-1} \sum_{\la \in \SP}
f(\la) q^{|\la|}.$$
Since $\langle I \rangle^{\SP}_q =1$, we may think of $\langle f
\rangle^{\SP}_q$ as the $q$-expectation value of $f$. We regard
$t^{\la_k}$ $(k \ge 1)$ as a function on $\SP$ (which by our
convention takes value $1$ on $\la$ with $\ell (\la)<k$).

\begin{lemma}   \label{lem:one}
We have
$$\sum_{\la\in \SP} t^{\la_1} q^{|\la|}
 =1 +\sum_{n=1}^\infty (1+q)(1+q^2) \cdots (1+q^{n-1}) q^{n} t^{n}.
 $$
\end{lemma}

\begin{proof}
Follows from the fact that the $n$-th term on the right-hand side
 of the identity equals the sum of $t^{\la_1} q^{|\la|}$ for $ \la_1
 =n$.
\end{proof}

The following lemma is a key step in our approach.

\begin{lemma}   \label{lem:general}
For $k \ge 1$, we have
$$\sum_{\la\in \SP} t^{\la_{k+1}} q^{|\la|}
 =\frac{q^{k(k+1)/2}}{(1-q)(1-q^2) \cdots (1-q^k)}
 \cdot \sum_{\la\in \SP} (q^k t)^{\la_1} q^{|\la|}.
 $$
\end{lemma}

\begin{proof}

By writing $\la =(a_1, \ldots, a_k, \mu)$ with $\mu =(\mu_1,
\mu_2,\ldots) \in \SP$ and $a_1   > \ldots > a_k >\mu_1$, we
obtain that
 \begin{eqnarray*}
 \sum_{\la\in \SP} t^{\la_{k+1}} q^{|\la|}
 &=& \sum_{a_1   > \ldots > a_k >\mu_1}
 \sum_{\mu\in \SP} t^{\mu_1} q^{a_1 +\ldots +a_k+|\mu|}.
 \end{eqnarray*}

Note also that (for a fixed $\mu_1$)
 \begin{eqnarray*}
  \sum_{a_1   > \ldots > a_k >\mu_1} q^{a_1 +\ldots +a_k}
   =\sum_{a_1   > \ldots > a_k >0} q^{a_1 +\ldots +a_k} q^{k
   \mu_1}.
 \end{eqnarray*}
 Thus the proof of the Lemma reduces to the following identity
 $$\sum_{a_1   > \ldots > a_k >0} q^{a_1 +\ldots +a_k}
 =\frac{q^{k(k+1)/2}}{(1-q)(1-q^2) \cdots (1-q^k)}.$$
This identity can be proved by a direct computation (as done in an
earlier version), or as the referee suggests, it is equivalent to
the simple fact that a strict partition of length $k$ can be
obtained from a partition of length $\leq k$ by adding row by row
$(k,k-1, \cdots, 1)$.
%
\end{proof}

\begin{theorem}  \label{th:firstform}
 We have
 \begin{eqnarray*}
 \left \langle \sum_{k \ge 1} t^{\la_k} \right \rangle^{\SP}_q
 = 1+ \sum_{n=1}^\infty \frac{q^n t^n}{1+ q^n}.
 \end{eqnarray*}
\end{theorem}

\begin{proof}

We recall a well-known identity of Euler (cf. e.g. Corollary~2.2,
\cite{An}):
\begin{eqnarray*}
 \sum_{k = 0}^\infty \frac{q^{k(k-1)/2} z^k}{(1-q)(1-q^2) \cdots (1-q^k)}
 = \prod_{r =0}^\infty (1 +q^{r} z).
\end{eqnarray*}
Here and below it is understood that the term for $k=0$ equals
$1$. We shall use an equivalent form of the Euler identity by
setting $z$ above to be $qz$:
\begin{eqnarray}  \label{eq:euler}
 \sum_{k = 0}^\infty \frac{q^{k(k+1)/2} z^k}{(1-q)(1-q^2) \cdots (1-q^k)}
 = \prod_{r =0}^\infty (1 +q^{r+1} z).
\end{eqnarray}

By Lemmas \ref{lem:one} and \ref{lem:general}, we have
 \begin{eqnarray}  \label{eq:formal}
 & & \sum_{\la\in \SP} \left ( \sum_{k \ge 0} t^{\la_{k+1}} \right )
   q^{|\la|} \nonumber \\
 &=& \sum_{k = 0}^\infty \frac{q^{k(k+1)/2}}{(1-q)(1-q^2) \ldots (1-q^k)}
 \cdot \sum_{\la\in \SP} (q^k t)^{\la_1} q^{|\la|}  \nonumber  \\
 &=& \sum_{k = 0}^\infty \frac{q^{k(k+1)/2}}{(1-q)(1-q^2) \cdots (1-q^k)}
   \nonumber \\
 && + \sum_{k = 0}^\infty \sum_{n=1}^\infty \frac{q^{k(k+1)/2}}{(1-q)(1-q^2) \cdots (1-q^k)}
 \cdot  (1+q)(1+q^2) \cdots (1+q^{n-1}) q^{n} (q^kt)^{n}  \nonumber \\
 &=& \sum_{k = 0}^\infty \frac{q^{k(k+1)/2}}{(1-q)(1-q^2) \cdots (1-q^k)}
   \nonumber \\
 & & + \sum_{n=1}^\infty (1+q)(1+q^2) \cdots (1+q^{n-1}) (qt)^{n}
 \sum_{k = 0}^\infty \frac{q^{k(k+1)/2}
 (q^{n})^k}{(1-q)(1-q^2) \cdots (1-q^k)}.
 \end{eqnarray}

Applying (\ref{eq:euler}) to the r.h.s. of (\ref{eq:formal}) with
$z=1$ and $z=q^n$, we have
 \begin{eqnarray*}
  && \sum_{\la\in \SP}
  \left(\sum_{k \ge 0} t^{\la_{k+1}}  \right) q^{|\la|}  \\
  &=& \prod_{r =0}^\infty (1 +q^{r+1})
   + \sum_{n=1}^\infty (1+q)(1+q^2) \cdots (1+q^{n-1}) (qt)^{n} \cdot
  \prod_{r =0}^\infty (1 +q^{n+r+1}) \\
  &=& (-q;q)_\infty + (-q;q)_\infty \sum_{n=1}^\infty \frac{q^n t^n}{1+ q^n}.
  \end{eqnarray*}
  This finishes the proof.
\end{proof}

\begin{remark}
 The same type of argument can be used to establish
 \begin{eqnarray}  \label{eq:super}
  && \sum_{\la\in \SP}
  \left(\sum_{k \ge 0} t^{\la_{k+1}} \right) q^{|\la|} z^{\ell(\la)}
  = (-qz;q)_\infty
   \left (1 + \sum_{n=1}^\infty \frac{q^n t^n z}{1+ q^nz} \right).
  \end{eqnarray}
 When $z=1$, it specializes to Theorem~\ref{th:firstform}. Another
 distinguished specialization is setting $z=-1$.

Our argument can also be adapted to give an alternative argument
for Theorem~6.5 of \cite{BO}, which is the counterpart of our
Theorem~\ref{th:firstform} above. The argument therein does not
seem to apply directly to our case, as the conjugation symmetry of
partitions used there is not available for strict partitions.
\end{remark}

By (\ref{eq:bridge}), we have
\begin{eqnarray}   \label{eq:relation}
:R:(t_1, \ldots, t_n) =(-q;q)_\infty \left \langle \prod_{i =1}^n
\sum_{k \ge 1} (t_i^{\la_k} -t_i^{-\la_k}) \right \rangle_q.
\end{eqnarray}
Therefore,
\begin{eqnarray*}
 R (t_1, \ldots, t_i, \ldots, t_n) &=& - R(t_1, \ldots, t_i^{-1}, \ldots t_n)  \\
 :R:(t_1, \ldots, t_i, \ldots, t_n) &=& - :R:(t_1, \ldots, t_i^{-1}, \ldots t_n).
\end{eqnarray*}

Recall the Jacobi $\theta$-functions:
$$\theta_{j,1} (q,t) =\sum_{n \in \frac{j}2 +\Z} q^{n^2} t^n,\quad
j=0,1.$$

Define
\begin{eqnarray}   \label{eq:thetaratio}
 \mathbb B (q,t) =   \frac{\theta_{1,1} (q, -t)}{\theta_{0,1} (q, -t)}
  = \left(- q^{-1/2}t \right)^{-\hf} \frac{(t; q^2)_\infty
(q^2t^{-1}; q^2)_\infty}{(qt; q^2)_\infty (qt^{-1} ; q^2)_\infty}.
\end{eqnarray}
The second equality above can be easily derived by using twice the
celebrated Jacobi triple product identity.

\begin{theorem}  \label{th:secondform}
 For $|q|<1, \; |q| < |t| < |q|^{-1}$,
\begin{enumerate}
\item the $1$-point correlation function $:R: (t)$ is given by
\begin{eqnarray*}
 :R: (t)
  &= & (-q;q)_\infty \sum_{n=1}^\infty \frac{q^n (t^n -t^{-n})}{1+ q^n} \\
  &= & (-q;q)_\infty \sum_{r =0}^\infty
  \left(
  \frac{(-1)^r q^{r+1} t}{1 -q^{r+1} t}
  -\frac{(-1)^r q^{r+1} t^{-1}}{1 -q^{r+1} t^{-1}}
 \right ).
 \end{eqnarray*}

 \item the one-point function
 $R(t)$ is given by
 \begin{eqnarray*}
 R(t) = (-q;q)_\infty \cdot t \frac{d}{dt} \ln \mathbb B (q, t).
 \end{eqnarray*}
 \end{enumerate}
\end{theorem}

\begin{proof}
Part (1) follows from (\ref{eq:relation}),
Theorem~\ref{th:firstform}, and the following simple identity:
 \begin{eqnarray*}
  \sum_{n=1}^\infty \frac{q^n t^n}{1+ q^n}
  =\sum_{r =0}^\infty \frac{(-1)^r q^{r+1} t}{1 -q^{r+1} t}.
  \end{eqnarray*}

 By definition of $R(t)$, we have
 \begin{eqnarray*}
 R(t)= :R:(t)
+\frac{t^{1/2}+t^{-1/2}}{2(t^{1/2}-t^{-1/2})} \cdot (-q;q)_\infty.
 \end{eqnarray*}

Now part (2) is a consequence of part (1) and the following simple
identities (and their counterparts with $t$ replaced by $t^{-1}$):
 \begin{eqnarray*}
  t \frac{d}{dt} \ln \left(t^{-1/2} (1-t) \right)
  &=&  \frac{t^{1/2}+t^{-1/2}}{2(t^{1/2}-t^{-1/2})} \\
   t \frac{d}{dt} \ln \left( (q^2t;q^2)_\infty (qt;q^2)_\infty^{-1} \right)
  &=& \sum_{r =0}^\infty \frac{(-1)^r q^{r+1} t}{1 -q^{r+1} t}.
  \end{eqnarray*}
  This finishes the proof.
\end{proof}

It is straightforward to check that
\begin{eqnarray}  \label{eq:diff}
 \mathbb B(q, qt) = - \mathbb B(q,t)^{-1}
 \end{eqnarray}
and thus by Theorem~\ref{th:secondform} that
$$R(qt) =- R(t).$$

\begin{remark}
The modular transformation properties of  $\mathbb B(q,t)$ (by
letting $q=e^{2\pi i \tau}$ and $t =e^{2\pi i z}$) and thus of
$R(t)$ (which actually depends on $q$ as well) can be derived from
the well-known modular properties of the $\theta$-functions and $
(-q;q)_\infty$.
\end{remark}

\section{Difference equations for the correlation functions}
\label{sec:property}

In this section we derive a $q$-difference equation satisfied by
the correlation function $R(t_1, \ldots, t_n).$

\begin{theorem} \label{th:difference}
 The function $R(t_1,\ldots, t_n)$ satisfies the
 following $q$-difference equation:
\begin{eqnarray*}
  && R(qt_1, t_2, \ldots, t_n) \\
 &=&
  \sum_{s=0}^{n-1} \sum_{1<i_1 <\cdots <i_s \leq n}
 \sum_{\varepsilon_{i_a}={\pm 1}} (-1)^{1+s+\#\varepsilon}
 R(t_1 t_{i_1}^{\varepsilon_{i_1}} \cdots
 t_{i_s}^{\varepsilon_{i_s}},
 \ldots, \widehat{t}_{i_1},\ldots, \widehat{t}_{i_s}, \ldots)
\end{eqnarray*}
where $\#\varepsilon$ denotes the number of $-1$'s among
$\varepsilon_{i_1}, \ldots, \varepsilon_{i_s}$, and
$\widehat{t}_{i}$ means that the very term is removed.
\end{theorem}

\begin{proof}

By (\ref{eq:action}) and (\ref{eq:keyop}), we have
\begin{eqnarray} \label{eq:commutator}
 [ \mathcal B_0(t), \phi_k ] = -(t^k -t^{-k}) \phi_k, \qquad k \in\Z.
\end{eqnarray}

Therefore,
\begin{eqnarray} \label{eq:perm}
 \bo (t_2) \cdots \bo (t_n)\phi_k
 = \sum_{S \subset \{2, \ldots, n\}} (-1)^{|S|}
 \prod_{i \in S} (t_i^k -t_i^{-k})\cdot
 \phi_k \prod_{i \not\in S} \bo (t_i).
\end{eqnarray}
This further implies that
\begin{eqnarray}  \label{eq:next}
 & &\Tr \mid_{\fock_0} q^{L_0} \phi_{-k} \bo (t_2) \cdots \bo (t_n)\phi_k \nonumber \\
 &=& \sum_{S \subset \{2, \ldots, n\}} (-1)^{|S|}
 \prod_{i \in S} (t_i^k -t_i^{-k}) \cdot \Tr \mid_{\fock_0} q^{L_0} \phi_{-k}\phi_k
 \prod_{i \not\in S} \bo (t_i).
\end{eqnarray}

On the other hand, by using $\phi_k q^{L_0}=q^{k}q^{L_0} \phi_k$
and the cyclic property of a trace, we have
\begin{eqnarray}  \label{eq:inv}
  \Tr \mid_{\fock_0} q^{L_0} \phi_{-k}
  \bo (t_2) \cdots \bo (t_n)\phi_k
 &=& \Tr \mid_{\fock_0} \phi_k q^{L_0}
 \phi_{-k} \bo (t_2) \cdots \bo (t_n) \nonumber \\
 &=& q^k \Tr \mid_{\fock_0}  q^{L_0}
  \phi_k\phi_{-k} \bo (t_2) \cdots \bo (t_n).
\end{eqnarray}

Define
$$\widetilde{\bo} (t) =\sum_{k \in \Z} t^k \phi_k \phi_{-k}.$$

By multiplying both sides of (\ref{eq:next}) with $t_1^k$ and then
summing over $k$, we obtain with the help of (\ref{eq:inv}) that
\begin{eqnarray*}
 & & \Tr \mid_{\fock_0} q^{L_0} \widetilde{\bo} (qt_1)
  \bo (t_2) \cdots \bo (t_n)  \\
 &=& \sum_{s=0}^{n-1} \sum_{1<i_1 <\cdots <i_s \leq n}
 \sum_{\varepsilon_{i_a}={\pm 1}} (-1)^{s+\#\varepsilon}
 R(t_1 t_{i_1}^{\varepsilon_{i_1}} \cdots
 t_{i_s}^{\varepsilon_{i_s}},
 \ldots, \widehat{t}_{i_1},\ldots, \widehat{t}_{i_s}, \ldots).
\end{eqnarray*}

Noting that
\begin{eqnarray}  \label{eq:correction}
 \bo (t) &=& :\bo(t):
 +\frac{t^{1/2}+t^{-1/2}}{2(t^{1/2}-t^{-1/2})} \text{I},
 \qquad |t| >1, \\
 \widetilde{\bo} (t) &=& - :\bo(t):
 -\frac{t^{1/2}+t^{-1/2}}{2(t^{1/2}-t^{-1/2})} \text{I},
 \qquad |t| <1, \nonumber
\end{eqnarray}
we have
\begin{eqnarray*}
\Tr \mid_{\fock_0} q^{L_0} \widetilde{\bo} (qt_1)
  \bo (t_2) \cdots \bo (t_n)
  = -R(qt_1, t_2, \ldots, t_n).
\end{eqnarray*}

This finishes the proof.
\end{proof}

\begin{remark}
Our formulation indicates clearly the nature of poles when
regarding $R(t_1,\ldots, t_n)$ as a meromorphic function. For
example, (\ref{eq:correction}) implies that
$$R(t_1, t_2, \ldots, t_n)
  =\frac{t_1 +1}{2(t_1 -1)}
  R(t_2, \ldots, t_n) + \text{\rm regular terms on } (t_1 =1). $$

As in \cite{BO, Ok}, the singularities of $R(t_1,\ldots, t_n)$ and
the $q$-difference functions determine $R(t_1,\ldots, t_n)$.
Indeed this is the key strategy used in {\em loc. cit} to
establish a compact closed formula for their $n$-point correlation
functions. It is expected the singularities of $R(t_1,\ldots,
t_n)$ and the $q$-difference equation can be very useful
eventually in finding a closed formula for $R(t_1,\ldots, t_n)$.
Similar remarks apply to the function $S(t_1,\ldots, t_n)$ studied
below.
\end{remark}




\begin{remark}  \label{rem:super}
 In light of (\ref{eq:super}), we can introduce a
 generalization of the $n$-point function $R(t_1,\ldots,t_n)$ by
 considering
 \begin{eqnarray*}
 R(t_1,\ldots,t_n; z) =\Tr\mid_{\fock_0} q^{L_0} z^{\alpha_0}
  \mathcal B_0(t_1)\cdots  \mathcal B_0(t_n),
 \end{eqnarray*}
where $\alpha_0 \stackrel{\text{def}}{=} \sum_{k >0} \phi_{-k}
\phi_k$. The specialization to $z=-1$, denoted by
$R^-(t_1,\ldots,t_n)$, might be regarded as a {\em super}
$n$-point function, which is no less natural to be considered than
the original $n$-point functions. It can be shown in a way similar
to Theorem~\ref{th:difference} to satisfy a $q$-difference
equation:
\begin{eqnarray*}
    R^-(qt_1, t_2, \ldots, t_n) = -R^-(t_2, \ldots, t_n) +
    \qquad \qquad \qquad \qquad \qquad \qquad \qquad \qquad \\
  \sum_{s=0}^{n-1} \sum_{1<i_1 <\cdots <i_s \leq n}
 \sum_{\varepsilon_{i_a}={\pm 1}} (-1)^{s+\#\varepsilon}
 R^-(t_1 t_{i_1}^{\varepsilon_{i_1}} \cdots
 t_{i_s}^{\varepsilon_{i_s}},
 \ldots, \widehat{t}_{i_1},\ldots, \widehat{t}_{i_s}, \ldots).
\end{eqnarray*}

Noting
 $$ [ \alpha_0, \phi_{\pm k} ] =\mp  \phi_{\pm k}, \quad k>0,$$
we see by (\ref{eq:super}) that (for $n=1$)
 \begin{eqnarray} \label{eq:modify}
 R(t; z) &=& \sum_{\la \in \SP}
\left (\sum_{k \ge 1} (t^{\la_k} -t^{-\la_k}) \right )
 q^{|\la|} z^{\ell (\la)}
 +(-qz;q)_\infty \frac{t^{1/2}+t^{-1/2}}{2(t^{1/2}-t^{-1/2})}
   \nonumber  \\
 &=& (-qz;q)_\infty \left (
 \sum_{n=1}^\infty \left(\frac{q^n t^n z}{1+ q^nz}
         - \frac{q^n t^n z}{1+ q^nz} \right)
 + \frac{t^{1/2}+t^{-1/2}}{2(t^{1/2}-t^{-1/2})}\right).
 \end{eqnarray}
In particular, by using (\ref{eq:modify}) we can show as in the
proof of Theorem~\ref{th:secondform} that
 $$R^-(t ) =  (q;q)_\infty \cdot t\frac{d}{dt}
 \ln \left (t^{-\hf}(t;q)_\infty (qt^{-1};q)_\infty \right ).
 $$

 Similar remarks apply to the function $S(t_1,\ldots,t_n)$ studied
 below.
\end{remark}

\section{Correlation functions on the twisted Fock space of Neveu-Schwarz type}
\label{sec:fock2}

Introduce a fermionic field of Neveu-Schwarz type

$$\varphi (z) =\sum_{n \in \hZ} \varphi_n z^{-n}$$
with the following commutation relations:
\begin{eqnarray*}
 [ \varphi_m, \varphi_n ]_+ = \delta_{m,-n}, \quad m, n \in \hZ.
\end{eqnarray*}

Denote by $\fockh$ the Fock space of $\varphi (z)$ with the
highest weight vector $\vac$ annihilated by $\varphi_n$, $n >0$. A
partition $\la =(\la_1, \la_2, \ldots \la_\ell)$ is called {\em
odd strict} if $\la_1 > \la_2 > \ldots > \la_\ell>0$ and all
$\la_i$ are odd integers. Introduce the notations
\begin{eqnarray*}
\varphi_\la &=& \varphi_{-\la_1/2} \varphi_{-\la_2/2} \ldots
\varphi_{-\la_\ell/2} \vac.
\end{eqnarray*}
Then $\fockh$ has a linear basis given by $\varphi_\la$, where
$\la$ runs over the set $\OSP$ of all odd strict partitions
(including the empty partition). The following grading operator
$L_0$ defines a natural $\hf \Z_+$-grading on $\fockh$:
$$L_0 (\varphi_\la)
= ({|\la|}/2)  \cdot \varphi_\la, \quad \la \in \OSP.$$

A subalgebra $\dinf$ of $\widehat{gl}_\infty$ acts on $\fockh$ as
follows \cite{DJKM} (also cf. \cite{KWY}). The algebra $\dinf$ is
spanned by $E_{i,j} -E_{1-j,1-i}$, where $i,j\in\Z$, and a central
element $C$, and the action is given concisely in terms of the
generating field
$$  \sum_{i,j \in \Z} (E_{i,j} -E_{1-j,1-i}) z^{i-\hf} w^{-j+\hf}
 = :\varphi (z) \varphi (w):
$$
and $C=I$. The Fock space $\fockh$ decomposes into a direct sum of
two irreducible representations of $\dinf$:
$$\fockh =\fockh_0 \bigoplus \fockh_1.$$
Here $\fockh_0$ (resp. $\fockh_1$) are spanned by $\varphi_{\la}$
with $\ell (\la)$ even (resp. odd).

Introduce the following operator $\mathcal D_0(t)$ in $\dinf$:
\begin{eqnarray*}
  \mathcal D_0 (t) = \sum_{i \in \hZ_+}   \left( t^{i-\hf} -t^{\hf-i}
   \right)
  (E_{i,i} -E_{1-i,1-i}) + \frac{t^{1/2}}{t - 1} I.
\end{eqnarray*}

The operator $\mathcal D_0(t)$ can be formally understood as
$$\mathcal D_0(t)= \sum_{k \in \hZ} t^{k} \varphi_{-k}\varphi_k$$
which is the zero-mode of $\varphi(tz) \varphi(z)$ (without normal
ordering). Indeed, for $|t|>1$, we have
$$\mathcal D_0 (t) = :\mathcal D_0 (t):
+ \left( \sum_{k\in \hZ, k<0} t^{k} \right )\text{I}
  =\; :\mathcal D_0 (t): +  \frac{t^\hf}{t-1} \; \text{I}$$
where
$$:\mathcal D_0 (t):
 = \sum_{k \in \hZ} t^{k} :\varphi_{-k}\varphi_k:.$$

We introduce the {\em $n$-point correlation functions}
\begin{eqnarray*}
 S(t_1, \ldots, t_n) &=&\Tr \mid_{\fockh} q^{L_0} \mathcal D_0
(t_1) \cdots \mathcal D_0 (t_n)
\end{eqnarray*}
and similarly define
\begin{eqnarray*}
 :S:(t_1, \ldots, t_n) &=& \Tr \mid_{\fockh} q^{L_0} :\mathcal D_0
(t_1): \cdots :\mathcal D_0 (t_n):.
\end{eqnarray*}
By definition, there is a simple identity of relating $S(t_1,
\ldots, t_n)$ to $:S:(t_1, \ldots, t_n)$. It is clear for $\la \in
\OSP$ that

\begin{eqnarray}  \label{eq:bridge2}
:\mathcal D_0 (t): (\varphi_\la)  = \sum_{k \ge 1} (t^{\la_k/2}
-t^{-\la_k/2}) \cdot \varphi_\la.
\end{eqnarray}
\section{The correlation function via odd strict partitions}
\label{sec:oddstr}

In this section, we develop the counterpart of
Sections~\ref{sec:strict} and \ref{sec:property}. We omit the
proofs which are similar to the earlier cases.

Note that
\begin{eqnarray*}
 \sum_{\la \in \OSP} q^{|\la|/2} =(-q^{1/2};q)_\infty.
\end{eqnarray*}
Given a function $f$ on the set $\OSP$ of odd strict partitions,
we denote by
$$\langle f \rangle^\OSP_q =(-q^{1/2};q)_\infty^{-1} \sum_{\la \in \OSP}
f(\la) q^{|\la|/2}.$$
Since $\langle I \rangle^\OSP_q =1$, we may think of $\langle f
\rangle^\OSP_q$ the $q$-expectation value of $f$. We regard
$t^{\la_k}$ $(k \ge 1)$ as a function on $\OSP$ which in
particular takes value $1$ on $\la$ with $\ell (\la)<k$.

\begin{lemma}   \label{lem:one2}
We have
$$\sum_{\la\in \OSP} t^{\la_1/2} q^{|\la|/2}
 =1 +\sum_{n=1}^\infty (1+q^{1/2})(1+q^{3/2}) \cdots (1+q^{n-1/2}) q^{n-1/2} t^{n-1/2}.
 $$
\end{lemma}

\begin{proof}
 Follows from the fact that the $n$-th term on the right-hand side
 of the identity equals the sum of $t^{\la_1/2} q^{|\la|/2}$ for $ \la_1
 =2n-1$.
\end{proof}

\begin{lemma}
For $k \ge 1$, we have
$$\sum_{\la\in \OSP} t^{\la_{k+1}/2} q^{|\la|/2}
 =\frac{q^{k(k+1)/2}}{(1-q)(1-q^2) \cdots (1-q^k)}
 \cdot \sum_{\la\in \OSP} (q^{k} t)^{\la_1/2} q^{|\la|/2}.
 $$
\end{lemma}

\begin{proof}
Follows from the same type of argument as in the proof of
Lemma~\ref{lem:general}.
\end{proof}

\begin{theorem}  \label{th:firstform2}
 We have
 \begin{eqnarray*}
 \left \langle \sum_{k \ge 1} t^{\la_k/2} \right \rangle^\OSP_q
  &=& 1+ \sum_{n=1}^\infty \frac{q^{n-1/2} t^{n-1/2}}{1+ q^{n-1/2}}
 \\
  \sum_{\la\in \OSP}
  \left(\sum_{k \ge 1} t^{\la_{k}/2} \right) q^{|\la|/2} z^{\ell(\la)}
  &=& (-q^{1/2}z)_\infty
   \left (1 + \sum_{n=1}^\infty \frac{q^{n-1/2}
   t^{n-1/2}z}{1+ q^{n-1/2}z} \right).
 \end{eqnarray*}
\end{theorem}

\begin{proof}
Follows from the same type of argument as in the proof of
Theorem~\ref{th:firstform}.
\end{proof}

By (\ref{eq:bridge2}), we have
\begin{eqnarray*}
:S:(t_1, \ldots, t_n) = (-q^{1/2};q)_\infty \left \langle \prod_{i
=1}^n \sum_{k \ge 1} (t_i^{\la_k} -t_i^{-\la_k}) \right
\rangle^\OSP_q.
\end{eqnarray*}

Recall $\mathbb B (q,t)$ was defined in (\ref{eq:thetaratio}).

\begin{theorem}  \label{th:secondform2}
 For $|q|<1, \; |q| < |t| < |q|^{-1}$,

 \begin{enumerate}
 \item the $1$-point correlation function $:S: (t)$ is given by
 \begin{eqnarray*}
 :S: (t)
  &=& (-q^{1/2};q)_\infty\cdot \sum_{n=1}^\infty \frac{q^{n-1/2}
  (t^{n-1/2} -t^{-n+1/2})}{1+ q^{n-1/2}} \\
  &=& (-q^{1/2};q)_\infty\cdot \sum_{r=0}^\infty
   \left (\frac{(-1)^r (q^{r+1} t)^{\hf}}{1 -q^{r+1} t}
  -\frac{(-1)^r (q^{r+1} t^{-1})^{\hf}}{1 -q^{r+1}  t^{-1}}
   \right )
 \end{eqnarray*}

 \item the $1$-point correlation function $S(t)$ is given by
 \begin{eqnarray*}
%
 S(t) = (-q^{1/2};q)_\infty \cdot t \frac{d}{dt} \left(
    \ln  \frac{\mathbb B (q^{1/2}, t^{1/2})}{\mathbb B (q^{1/2}, -t^{1/2})}
    \right).
 \end{eqnarray*}
 \end{enumerate}
\end{theorem}

\begin{proof}
Part (1) follows from Theorem~\ref{th:firstform2} and the
following simple identity:
 \begin{eqnarray*}
  \sum_{n=1}^\infty \frac{q^{n-1/2} t^{n-1/2} }{1+ q^{n-1/2} }
  &=& \sum_{r =0}^\infty \frac{(-1)^r (q^{r+1} t)^{\hf}}{1 -q^{r+1}
  t}.
 \end{eqnarray*}.

By definition of $S(t)$, we have
 \begin{eqnarray*}
S(t)= :S:(t) +\frac{t^{\hf}}{t -1} \cdot (-q^{1/2};q)_\infty.
 \end{eqnarray*}

We can rewrite by the definition (\ref{eq:thetaratio}) of $\mathbb
B$ that
\begin{eqnarray} \label{eq:series}
 && \ln  {\mathbb B (q^{1/2}, t^{1/2})} -\ln{\mathbb B
 (q^{1/2}, -t^{1/2})}   \nonumber  \\
 &=& \frac{1-t^{1/2}}{1+t^{1/2}}
 \prod_{k \ge 1} (1-q^{k/2} t^{1/2})^{(-1)^k}
  (1 +q^{k/2} t^{1/2})^{(-1)^{k+1}} \cdot \nonumber   \\
 && \qquad \qquad \cdot \prod_{k \ge 1} (1-q^{k/2} t^{-1/2})^{(-1)^k}
  (1 +q^{k/2} t^{-1/2})^{(-1)^{k+1}}.
 \end{eqnarray}

Part (2) can now be derived from (\ref{eq:series}) and  part (1)
of the Theorem.
\end{proof}

It follows from (\ref{eq:diff}) and Theorem~\ref{th:secondform2}
that $S(qt) =- S(t).$
%
The modular transformation properties of $S(t)$ can be derived
from the well-known modular properties of the $\theta$-functions
and $(-q^{1/2};q)_\infty$.

%
%

The same proof as in Theorem~\ref{th:difference} also leads to the
following $q$-difference equation for $S$. It is interesting to
note that the two distinct functions $R(t_1,\ldots, t_n)$ and
$S(t_1,\ldots, t_n)$ satisfy the same difference equations but are
distinguished by their singularities.

\begin{theorem}
 The meromorphic function $S(t_1,\ldots, t_n)$ satisfies the
 following difference equation:
\begin{eqnarray*}
  && S(qt_1, t_2, \ldots, t_n) \\
 &=&
  \sum_{s=0}^{n-1} \sum_{1<i_1 <\cdots <i_s \leq n}
 \sum_{\varepsilon_{i_a}={\pm 1}} (-1)^{1+s+\#\varepsilon}
 S(t_1 t_{i_1}^{\varepsilon_{i_1}} \cdots
 t_{i_s}^{\varepsilon_{i_s}},
 \ldots, \widehat{t}_{i_1},\ldots, \widehat{t}_{i_s}, \ldots)
\end{eqnarray*}
where $\#\varepsilon$ denotes the number of $-1$'s among
$\varepsilon_{i_1}, \ldots, \varepsilon_{i_s}$.
\end{theorem}
\section{Discussions} \label{sec:discuss}

In this paper we have studied certain distinguished operators on
twisted Fock spaces and the corresponding correlation functions.
These functions satisfy certain $q$-difference equations and they
can be interpreted as $q$-expectation values of functions on the
set of (odd) strict partitions. We found closed formulas for the
$1$-point functions in terms of Jacobi $\theta$-functions etc,
which implies nice modular transformation properties.

The Lie algebra $\winfty$ of differential operators on the circle
has two distinguished subalgebras $\widehat{\mathcal D}^\pm$
induced from anti-involutions \cite{KWY}. There exists canonical
Lie algebra homomorphisms from $\widehat{\mathcal D}^\pm$ to the
Lie algebras $\binf$, $\dinf$. In this way, the Fock spaces
$\fock$ and $\fockh$ can be regarded as representations of these
Lie algebras $\widehat{\mathcal D}^\pm$, and they decompose into a
direct sum of two irreducibles. The study of correlation functions
in this paper can thus be interpreted as character formulas of
these $\widehat{\mathcal D}^\pm$-modules. This is similar to the
case of \cite{BO}, where the correlation functions of the infinite
wedge representations were also interpreted as characters of
$\winfty$.

In \cite{FKRW} (also cf. \cite{AFMO} and the references therein),
irreducible quasi-finite modules of $\winfty$ of higher levels
using Fock space were studied in detail, and a similar study has
been carried out \cite{KWY} for irreducible quasi-finite modules
of $\widehat{\mathcal D}^\pm$ of higher levels using untwisted and
twisted Fock spaces. It will be interesting to study the
correlation functions on these modules of higher levels.

In is an interesting question to find closed formulas for the
correlation functions $R(t_1, \ldots, t_n)$ and $S(t_1, \ldots,
t_n)$ when $n>1$. The $q$-difference equations and singularities
of the functions $R$ and $S$ are expected to play an important
role. While the vertex operators are not essentially used in the
paper, the correlation functions here can be related to the
correlation functions in vertex algebras (cf. \cite{Zhu, DLM,
Mi}), and this connection could be useful in understanding these
functions further. We refer to Milas \cite{Mi} for very
interesting results in this direction.

\end{document}